\documentclass[a4paper,12pt,reqno]{amsart}
\usepackage{amsmath}
\usepackage{amssymb}
\usepackage{amsthm}

\usepackage{amscd}
\usepackage{amsfonts}
\usepackage{setspace}
\usepackage{version}

\usepackage{hyperref}

\title[Moments of random multiplicative functions review]{Moments of random multiplicative functions, III: A short review}
\author{Adam J Harper}
\address{Mathematics Institute, Zeeman Building, University of Warwick, Coventry CV4 7AL, England}
\email{A.Harper@warwick.ac.uk}
\date{15th October 2024}
\thanks{This research was funded in part by the Engineering and Physical Sciences Research Council of the United Kingdom [grant EP/V055755/1]. For the purpose of open access, the author has applied a Creative Commons Attribution (CC-BY) licence to any Author Accepted Manuscript version arising from this submission.}

\numberwithin{equation}{section}
\theoremstyle{plain}

\onehalfspace
\newcommand{\N}{\mathbb{N}}
\newcommand{\R}{\mathbb{R}}
\newcommand{\E}{\mathbb{E}}
\newcommand{\p}{\mathbb{P}}
\newcommand{\Z}{\mathbb{Z}}
\newcommand{\C}{\mathbb{C}}

\newtheorem{theorem}{Theorem}
\numberwithin{theorem}{section}

\numberwithin{proposition}{section}

\numberwithin{corollary}{section}

\begin{document}

\maketitle

\begin{abstract}
We give a short review of recent progress on determining the order of magnitude of moments $\E|\sum_{n \leq x} f(n)|^{2q}$ of random multiplicative functions, and of closely related issues.

We hope this can serve as a concise introduction to some of the ideas involved, for those who may not have too much background in the area.
\end{abstract}

\section{Introduction}
{\em Preliminary remark on notation.} We write $f(x) = O(g(x))$ and $f(x) \ll g(x)$, both meaning (as usual in analytic number theory, but perhaps not elsewhere) that there exists $C$ such that $|f(x)| \leq Cg(x)$, for all $x$. Sometimes this notation will be adorned with a subscript parameter (e.g. $O_{\epsilon}(\cdot)$ and $\ll_{\delta}$), meaning that the implied constant $C$ may depend on that parameter. We write $f(x) \asymp g(x)$ to mean that $g(x) \ll f(x) \ll g(x)$, in other words that $c g(x) \leq |f(x)| \leq C g(x)$ for some $c,C$, for all $x$.

\vspace{12pt}
Let $(f(p))_{p \; \text{prime}}$ be a sequence of independent Steinhaus random variables, i.e. independent random variables distributed uniformly on the complex unit circle $\{|z|=1\}$. We define a {\em Steinhaus random multiplicative function} $f : \N \rightarrow \C$, by setting $f(n) := \prod_{p^{a} || n} f(p)^{a}$ for all natural numbers $n$ (where $p^a || n$ means that $p^a$ is the highest power of the prime $p$ that divides $n$, so $n = \prod_{p^{a} || n} p^{a}$). Thus $f$ is a random function taking values in the complex unit circle, that is totally multiplicative, i.e. satisfies $f(nm) = f(n)f(m)$ for all $n,m$.

For simplicity, in this survey we shall confine our attention to these Steinhaus random multiplicative functions, and not discuss other models (e.g. the Rademacher or extended Rademacher models).

\vspace{12pt}
Random multiplicative functions sit at the intersection of number theory, probability, and analysis. Thus Steinhaus random multiplicative functions provide a heuristic model for randomly chosen Dirichlet characters $\chi(n)$ or ``continuous characters'' $n \mapsto n^{it}$: see e.g. the papers of Granville and Soundararajan~\cite{gransoundlcs} and Lamzouri~\cite{lamzouri2dzeta}. In some circumstances, they also serve as tools for proving rigorous results about such number theoretic objects (see e.g. the author's paper~\cite{harpertypicalchar}). From a probabilistic point of view, the values of a Steinhaus random multiplicative function are a naturally arising sequence of {\em dependent} random variables: notice e.g. that $f(6)=f(2)f(3)$, so the triple of values $f(2), f(3), f(6)$ are clearly not all independent of one another. Then one wishes to understand how this dependence influences the behaviour, compared with the classical probabilistic setting of sequences of independent random variables.

\vspace{12pt}
In this survey we shall describe these efforts from the specific perspective of the (absolute) {\em power moments} $\E|\sum_{n \leq x} f(n)|^{2q}$, where $q \geq 0$ is real. The only cases that are really easy to handle are the trivial case $q=0$, and the second moment case $q=1$. For we can observe that
$$ \E f(n) \overline{f(m)} = \E \left(\prod_{p^{a} || n} f(p)^{a} \right) \left(\prod_{p^{a} || m} f(p)^{-a} \right) = \E \prod_{p| mn} f(p)^{a(n,p) - a(m,p)} = \prod_{p| mn} \E f(p)^{a(n,p) - a(m,p)} , $$
where $a(n,p), a(m,p)$ are the exponents of $p$ in the unique prime factorisations of $n,m$, respectively. Since $f(p)^{a(n,p) - a(m,p)}$ is uniformly distributed on the unit circle (and in particular has mean zero) except when $a(n,p) = a(m,p)$, it follows that $\E f(n) \overline{f(m)} = \textbf{1}_{a(n,p)=a(m,p) \; \forall \, p} = \textbf{1}_{n=m}$, where $\textbf{1}$ is the indicator function\footnote{In other words, the sequence of random variables $(f(n))_{n \in \N}$ are {\em orthogonal}.}. Thus
$$ \E|\sum_{n \leq x} f(n)|^{2} = \E \sum_{n,m \leq x} f(n) \overline{f(m)} = \sum_{n,m \leq x} \textbf{1}_{n=m} = \lfloor x \rfloor . $$

Since the second moment has size $\approx x$, the most immediate simple conjecture (e.g. thinking of the moments of Gaussian random variables, or of sums of independent random variables) might be that the $2q$-th moment should have size $\asymp_{q} x^q$. But the true behaviour is far more subtle and interesting.

\begin{theorem}[Harper~\cite{harperrmflowmoments}, 2020]\label{thmstlow}
If $f(n)$ is a Steinhaus random multiplicative function, then uniformly for all large $x$ and all real $0 \leq q \leq 1$ (possibly depending on $x$) we have
$$ \E|\sum_{n \leq x} f(n)|^{2q} \asymp \left( \frac{x}{1 + (1-q)\sqrt{\log\log x}} \right)^q . $$
\end{theorem}

\begin{theorem}[Harper~\cite{harperrmfhigh}, 2019]\label{thmsthigh}
There exist a small absolute constant $c > 0$, and a large absolute constant $C > 0$, such that the following is true. If $f(n)$ is a Steinhaus random multiplicative function, then uniformly for all large $x$ and real $1 \leq q \leq \frac{c\log x}{\log\log x}$ we have
$$ e^{-q^{2}\log q - q^{2}\log\log(2q) - Cq^2} \leq \frac{\E |\sum_{n \leq x} f(n)|^{2q}}{x^{q} \log^{(q-1)^2}x} \leq e^{-q^{2}\log q - q^{2}\log\log(2q) + Cq^2} . $$
\end{theorem}

Theorem \ref{thmstlow} implies that $\E|\sum_{n \leq x} f(n)| \asymp \frac{\sqrt{x}}{(\log\log x)^{1/4}}$, which resolved a conjecture of Helson~\cite{helson} that the first absolute moment should be $o(\sqrt{x})$. For any positive $\lambda$, Markov's inequality and Theorem \ref{thmstlow} also immediately yield that
$$ \p(|\sum_{n \leq x} f(n)| > \lambda \frac{\sqrt{x}}{(\log\log x)^{1/4}}) \leq \frac{\E|\sum_{n \leq x} f(n)|}{\lambda \frac{\sqrt{x}}{(\log\log x)^{1/4}}} \ll \frac{1}{\lambda} . $$
So we may say that typically (e.g. with probability $\geq 0.99$) the sums $\sum_{n \leq x} f(n)$ are $\ll \frac{\sqrt{x}}{(\log\log x)^{1/4}} = o(\sqrt{x})$, enjoying {\em better than squareroot cancellation} (or {\em subdiffusivity}, in more probabilistic language), as opposed to the squareroot size suggested by the second moment. This is quite rare and unexpected in number theoretic settings. See \cite{harperrmflowmoments,harperrmfhigh} for more precise bounds on the large deviations of $\sum_{n \leq x} f(n)$. Theorem \ref{thmstlow} and H\"older's inequality also directly imply that $|\sum_{n \leq x} f(n)| \gg \frac{\sqrt{x}}{(\log\log x)^{1/4}}$ with positive probability, since for any small parameter $\eta > 0$ we have
\begin{eqnarray}
\frac{\sqrt{x}}{(\log\log x)^{1/4}} & \ll & \E|\sum_{n \leq x} f(n)| \leq \frac{\eta \sqrt{x}}{(\log\log x)^{1/4}} + \E\textbf{1}_{|\sum_{n \leq x} f(n)| > \eta \frac{\sqrt{x}}{(\log\log x)^{1/4}}}|\sum_{n \leq x} f(n)| \nonumber \\
& \leq & \frac{\eta \sqrt{x}}{(\log\log x)^{1/4}} + \p(|\sum_{n \leq x} f(n)| > \eta \frac{\sqrt{x}}{(\log\log x)^{1/4}})^{1/3} (\E|\sum_{n \leq x} f(n)|^{3/2})^{2/3} \nonumber \\
& \ll & \frac{\eta \sqrt{x}}{(\log\log x)^{1/4}} + \p(|\sum_{n \leq x} f(n)| > \eta \frac{\sqrt{x}}{(\log\log x)^{1/4}})^{1/3} \frac{\sqrt{x}}{(\log\log x)^{1/4}} . \nonumber
\end{eqnarray}
Provided $\eta$ is fixed sufficiently small, this forces $\p(|\sum_{n \leq x} f(n)| > \eta \frac{\sqrt{x}}{(\log\log x)^{1/4}}) \gg 1$.

Looking at Theorems \ref{thmstlow} and \ref{thmsthigh} qualitatively, we see that for fixed $q < 1$ the moments reflect the typical size $\frac{\sqrt{x}}{(\log\log x)^{1/4}}$ of $\sum_{n \leq x} f(n)$, whereas when $q=1$ the second moment (being of order $x$) does {\em not} reflect the typical size, instead (and perhaps unexpectedly) being dominated by somewhat unusual larger values of $\sum_{n \leq x} f(n)$. Note that in both theorems it is permissible to choose $q$ in a way that depends on $x$, so one can explore the nature of the transition when $q = 1+o(1)$. As $q$ becomes even larger, the moments are dominated by increasingly large and rare values of $\sum_{n \leq x} f(n)$, giving rise to the rapidly growing term $\log^{(q-1)^2}x$ in Theorem \ref{thmsthigh}. In particular, it is only the low moments ($q < 1$) that give access to the typical behaviour of $\sum_{n \leq x} f(n)$.

\vspace{12pt}
The cases of Theorem \ref{thmsthigh} where $q \in \N$ can be successfully attacked by expanding the $2q$-th power, using the orthogonality property $\E f(n_1) ... f(n_q) \overline{f(m_1) ... f(m_q)} = \E f(n_1 ... n_q) \overline{f(m_1 ... m_q)} = \textbf{1}_{n_1 ... n_q = m_1 ... m_q}$, and trying to bound or evaluate the divisor type sum that remains. See Harper, Nikeghbali and Radziwi{\l}{\l}~\cite{hnr}, and Heap and Lindqvist~\cite{heaplindqvist}, as well as unpublished work of Granville and Soundararajan. Aside from this, no sharp bounds were known in any case of Theorems \ref{thmstlow} and \ref{thmsthigh} prior to the work \cite{harperrmflowmoments,harperrmfhigh} of the author, and it is the approach and ideas from those papers (with a few later refinements and elaborations) that we shall try to explain in the following sections. See the introductions to \cite{harperrmflowmoments,harperrmfhigh} for further references to previously known, non-sharp moment bounds (lower bounds for the low moments, upper and lower bounds for high moments).

In the final section, we also provide a small selection of further reading on related topics.

\section{Reducing to random Euler products}\label{secep}
The classical number theoretic approach to studying sums of multiplicative functions entails introducing suitable multiplicative generating functions, like the Riemann zeta function or Dirichlet $L$-functions. For Theorems \ref{thmstlow} and \ref{thmsthigh}, a natural choice of generating function is
$$ F(s) := \sum_{\substack{n=1, \\ p|n \Rightarrow p \leq x}}^{\infty} \frac{f(n)}{n^s} = \prod_{\text{prime} \, p \leq x} (1 - \frac{f(p)}{p^{s}})^{-1} , $$
the {\em random Euler product} corresponding to $f(n)$. Since this is a finite product, it certainly converges whenever $\Re(s) > 0$ (so that $|\frac{f(p)}{p^{s}}| < 1$). An immediate appeal of $F(s)$ is that it takes the form of a product of {\em independent} factors $(1 - \frac{f(p)}{p^{s}})^{-1}$.

\vspace{12pt}
To connect $\sum_{n \leq x} f(n)$ with $F(s)$, the obvious route is to apply Perron's formula (multiplicative Fourier inversion), which would yield something like
$$ \sum_{n \leq x} f(n) \approx \frac{1}{2\pi i} \int_{\sigma - iT}^{\sigma + iT} F(s) \frac{x^s}{s} ds , \;\;\;\;\; \sigma > 0 , $$
for a suitably large parameter $T$. In principle it should then be possible to perfectly understand the distribution of $\sum_{n \leq x} f(n)$, and in particular to understand its moments, by perfectly understanding the (joint) distribution of $F(s)$ for various $s$, but in practice difficulties immediately arise. For example, if we want to estimate $\E|\sum_{n \leq x} f(n)|$, the only really obvious approach is to use the triangle inequality, obtaining that
$$ \E|\sum_{n \leq x} f(n)| \lesssim \frac{1}{2\pi} \int_{-T}^{T} \E|F(\sigma + it)| \frac{x^\sigma}{|\sigma + it|} dt . $$
Since we expect the left hand side to be around $\sqrt{x}$, the natural choice of abscissa is $\sigma = 1/2$, and we can imagine that $T \asymp 1$, say (it certainly cannot be smaller). But it is not hard to calculate that $\E|F(1/2+it)| \asymp \log^{1/4}x$, so we only obtain an upper bound $\E|\sum_{n \leq x} f(n)| \lesssim \sqrt{x} \log^{1/4}x$, which is significantly worse than the trivial Cauchy--Schwarz bound $\E|\sum_{n \leq x} f(n)| \leq \sqrt{\E|\sum_{n \leq x} f(n)|^2} \leq \sqrt{x}$.

The problem is that to understand the moments directly from Perron's formula, and in particular to have any hope of capturing the delicate double logarithmic saving in the low moments\footnote{When studying $\E|\sum_{n \leq x} f(n)|^{2q}$ with $q$ large enough, one {\em can} obtain sharp {\em upper bounds} by starting with Perron's formula and the triangle inequality, because for sufficiently high moments the main contribution to the Perron integral comes from just a few large values of $F(1/2+it)$ (at some random $t$), so the triangle inequality doesn't lose much. See the end of the introduction of Harper~\cite{harperrmfhigh}, and see Szab\'o's paper~\cite{szaboupper} for an implementation of similar ideas in the context of character sums.} in Theorem \ref{thmstlow}, one would need to understand the full value distribution (both modulus and argument) of all the $F(1/2+it)$, and the interaction of this with the phase $x^{it}$ that is destroyed by the triangle inequality. It seems very challenging to operate at such a level of precision, and so the papers~\cite{harperrmflowmoments,harperrmfhigh} adopt a less direct approach.

\vspace{12pt}
Before connecting with a random Euler product, we first work with $\sum_{n \leq x} f(n)$ ``by hand''. We shall outline the argument in a form roughly suitable for ultimately proving the upper bound part of Theorem \ref{thmstlow}, and then indicate the changes needed when working towards lower bounds in Theorem \ref{thmstlow}, and towards Theorem \ref{thmsthigh}.

Given any large parameter $P$, we say a number is $P$-{\em rough} if all of its prime factors are $> P$, and $P$-{\em smooth} if all of its prime factors are $\leq P$. Then using the multiplicativity of $f$, we may write
\begin{eqnarray}\label{rsdisplay}
\sum_{n \leq x} f(n) & = & \sum_{\substack{n \leq x, \\ n \; \text{has a prime factor} \; > P}} f(n) + \sum_{\substack{n \leq x, \\ n \; \text{is} \; P \; \text{smooth}}} f(n) \nonumber \\
& = & \sum_{\substack{P < m \leq x, \\ m \; \text{is} \; P \; \text{rough}}} f(m) \sum_{\substack{n \leq x/m, \\ n \; \text{is} \; P \; \text{smooth}}} f(n) + \sum_{\substack{n \leq x, \\ n \; \text{is} \; P \; \text{smooth}}} f(n) .
\end{eqnarray}
When $0 \leq q \leq 1$, H\"older's inequality and the orthogonality property of $f$ imply that $\E|\sum_{\substack{n \leq x, \\ n \; \text{is} \; P \; \text{smooth}}} f(n)|^{2q} \leq (\E|\sum_{\substack{n \leq x, \\ n \; \text{is} \; P \; \text{smooth}}} f(n)|^{2})^{q} = (\#\{n \leq x : n \; \text{is} \; P \; \text{smooth}\})^{q}$, and it is not too hard to show that $\#\{n \leq x : n \; \text{is} \; P \; \text{smooth}\} \ll x e^{-(\log x)/\log P}$ (and indeed much more precise estimates are known). So provided we choose $P$ with $\log P$ somewhat smaller than $\log x$ (e.g. $P = x^{1/\log\log x}$), the contribution to the moment from $\sum_{\substack{n \leq x, \\ n \; \text{is} \; P \; \text{smooth}}} f(n)$ will be negligible.

To work with the first sums in \eqref{rsdisplay}, we use one of the key general techniques in the study of random multiplicative functions, namely {\em conditioning}. If we let $\E^P$ denote expectation conditional on the values $(f(p))_{p \leq P}$ (i.e. expectation with those values treated as fixed and the $(f(p))_{p > P}$ remaining random, so the conditional expectation of any quantity is a function of the values $(f(p))_{p \leq P}$), then the Tower Property of conditional expectation implies that
$$ \E|\sum_{\substack{P < m \leq x, \\ m \; \text{is} \; P \; \text{rough}}} f(m) \sum_{\substack{n \leq x/m, \\ n \; \text{is} \; P \; \text{smooth}}} f(n)|^{2q} = \E \E^P |\sum_{\substack{P < m \leq x, \\ m \; \text{is} \; P \; \text{rough}}} f(m) \sum_{\substack{n \leq x/m, \\ n \; \text{is} \; P \; \text{smooth}}} f(n)|^{2q} . $$
Note that in this case, the Tower Property is simply Fubini’s theorem, breaking up the multiple ``integration'' $\E$ into separate integrations corresponding to the $(f(p))_{p \leq P}$ (on the outside) and the $(f(p))_{p > P}$. Then applying H\"older's inequality {\em to the conditional expectation $\E^P$ only}, and subsequently applying orthogonality of $f$ to evaluate the conditional second moment that emerges, we see the above is
\begin{equation}\label{condhold}
\leq \E \Biggl( \E^P |\sum_{\substack{P < m \leq x, \\ m \; \text{is} \; P \; \text{rough}}} f(m) \sum_{\substack{n \leq x/m, \\ n \; \text{is} \; P \; \text{smooth}}} f(n)|^{2} \Biggr)^q = \E \Biggl( \sum_{\substack{P < m \leq x, \\ m \; \text{is} \; P \; \text{rough}}} |\sum_{\substack{n \leq x/m, \\ n \; \text{is} \; P \; \text{smooth}}} f(n)|^{2} \Biggr)^q .
\end{equation}

These manipulations have allowed us to efficiently pass from working with $\sum_{n \leq x} f(n)$ at a single point $x$, to working with a mean square average of $\sum_{\substack{n \leq x/m, \\ n \; \text{is} \; P \; \text{smooth}}} f(n)$ over many points $x/m$. Some fairly simple {\em sieve theory}, coupled with a smoothing argument where the sum over $m$ is broken into smaller pieces that can be well approximated (on average) by integrals, shows the above is
$$ \approx \E \Biggl( \frac{1}{\log P} \int_{P}^{x} |\sum_{\substack{n \leq x/t, \\ n \; \text{is} \; P \; \text{smooth}}} f(n)|^{2} dt \Biggr)^q . $$
Here $\frac{1}{\log P}$ is the approximate density of the $P$-rough numbers, as revealed by sieve theory bounds\footnote{For background on this, an unfamiliar reader may consult Chapter 3 of Montgomery and Vaughan~\cite{mv}, for example.}. It is important in this argument that $\sum_{\substack{n \leq x, \\ n \; \text{is} \; P \; \text{smooth}}} f(n)$, corresponding to $m=1$, was separated and removed in \eqref{rsdisplay}, since for this piece one could not smooth the one-term $m$ ``sum'' and would not pick up a $\frac{1}{\log P}$ density saving.

Substituting $z = x/t$ in the integral, we find the above expression is
$$ =  x^q \E \Biggl( \frac{1}{\log P} \int_{1}^{x/P} |\sum_{\substack{n \leq z, \\ n \; \text{is} \; P \; \text{smooth}}} f(n)|^{2} \frac{dz}{z^2} \Biggr)^q . $$
Now the quantity inside the expectation is simply an integral mean square average of $\sum_{\substack{n \leq z, \\ n \; \text{is} \; P \; \text{smooth}}} f(n)$, so rather than inefficiently applying Perron's formula (multiplicative Fourier inversion) and the triangle inequality to connect with a random Euler product, we can efficiently apply the {\em multiplicative version of Parseval's identity} (with no phases $x^{it}$ to be destroyed). One obtains\footnote{In the original papers  \cite{harperrmflowmoments,harperrmfhigh}, this argument was run in a bit more complicated way, with \eqref{rsdisplay} replaced by a subdivision of $\sum_{n \leq x} f(n)$ into multiple subsums according to various possible ranges for the largest prime factor of $n$. Whilst a careful subdivision does seem to be necessary in the high moments case~\cite{harperrmfhigh}, it is not in the low moments case, and the cleaner argument we outlined here (with just a single parameter $P$) can be implemented rigorously. In a character sum context, this is essentially what is done by Harper~\cite{harpertypicalchar}. For random multiplicative functions, it is done in the recent preprint of Gorodetsky and Wong~\cite{gorwong}. \label{fnsubdiv}} a bound
$$ \leq x^q \E \Biggl( \frac{1}{\log P} \int_{1}^{\infty} |\sum_{\substack{n \leq z, \\ n \; \text{is} \; P \; \text{smooth}}} f(n)|^{2} \frac{dz}{z^2} \Biggr)^q =  x^q \E \Biggl( \frac{1}{2\pi \log P} \int_{-\infty}^{\infty} \frac{|F_{P}(1/2+it)|^{2}}{|1/2 + it|^2} dt \Biggr)^q , $$
where $F_P(s) := \sum_{\substack{n=1, \\ p|n \Rightarrow p \leq P}}^{\infty} \frac{f(n)}{n^s} = \prod_{p \leq P} (1 - \frac{f(p)}{p^{s}})^{-1}$.

As a final simplification, note that because the joint distribution of $(f(p)p^{-it})_{p \; \text{prime}}$ is the same for any fixed $t \in \R$ (namely a sequence of independent Steinhaus random variables), it follows that for any given set $\mathcal{H} \subseteq \R$, the joint distribution of $(F_{P}(1/2+ih+it))_{h \in \mathcal{H}}$ is the same for all $t \in \R$. We call this property {\em translation invariance in law}, and it implies that $\E \left( \int_{n-1/2}^{n+1/2} \frac{|F_{P}(1/2+it)|^{2}}{|1/2 + it|^2} dt \right)^q$ is
$$ \ll \frac{1}{1+|n|^{2q}} \E \Biggl( \int_{n-1/2}^{n+1/2} |F_{P}(1/2+it)|^{2} dt \Biggr)^q = \frac{1}{1+|n|^{2q}} \E \Biggl( \int_{-1/2}^{1/2} |F_{P}(1/2+it)|^{2} dt \Biggr)^q $$
for all $n \in \Z$. Thus provided $2/3 \leq q \leq 1$, say, (so that $\sum_{n \in \Z} \frac{1}{1+|n|^{2q}}$ converges), we have the simpler bound
$$ x^q \E \Biggl( \frac{1}{2\pi \log P} \int_{-\infty}^{\infty} \frac{|F_{P}(1/2+it)|^{2}}{|1/2 + it|^2} dt \Biggr)^q \ll x^q \E \Biggl( \frac{1}{\log P} \int_{-1/2}^{1/2} |F_{P}(1/2+it)|^{2} dt \Biggr)^q . $$
Thanks to H\"older's inequality, note that if we can prove Theorem \ref{thmstlow} on this range $2/3 \leq q \leq 1$ then we can immediately deduce it on the full range $0 \leq q \leq 1$ as well.

\vspace{12pt}
We end this section with a summary of the modifications required in these arguments when working towards lower bounds in Theorem \ref{thmstlow}, and towards Theorem \ref{thmsthigh}.
\begin{itemize}
\item For lower bounds on the low moments range $0 \leq q \leq 1$, the only substantial change comes at the beginning, where \eqref{condhold} currently goes in the wrong direction (giving an upper rather than lower bound). Instead, if we write $P(n)$ for the largest prime factor of $n$, then a fairly simple argument (essentially just the triangle inequality) shows that
$$ \E|\sum_{n \leq x} f(n)|^{2q} \gg \E|\sum_{\substack{n \leq x, \\ P(n) > x^{3/4}}} f(n)|^{2q} = \E|\sum_{x^{3/4} < p \leq x} f(p) \sum_{n \leq x/p} f(n)|^{2q} . $$
Using the Tower Property of conditional expectation, the right hand side here is $= \E \E^{x^{3/4}} |\sum_{x^{3/4} < p \leq x} f(p) \sum_{n \leq x/p} f(n)|^{2q}$. And when we condition on the values of $f$ on all primes $\leq x^{3/4}$, the innermost sums $\sum_{n \leq x/p} f(n)$ become fixed, so $\sum_{x^{3/4} < p \leq x} f(p) \sum_{n \leq x/p} f(n)$ is (under this conditioning) simply a {\em sum of independent random variables} $(f(p))_{x^{3/4} < p \leq x}$ multiplied by some coefficients. The moments of such classical sums are very well understood, for example Khintchine's inequality implies that
$$ \E^{x^{3/4}} |\sum_{x^{3/4} < p \leq x} f(p) \sum_{n \leq x/p} f(n)|^{2q} \gg \left( \sum_{x^{3/4} < p \leq x} |\sum_{n \leq x/p} f(n)|^2 \right)^q . $$
Thus $\E|\sum_{n \leq x} f(n)|^{2q} \gg \E \left( \sum_{x^{3/4} < p \leq x} |\sum_{n \leq x/p} f(n)|^2 \right)^q$, which provides a suitable lower bound analogue\footnote{One could replace $x^{3/4}$ by any value between $\sqrt{x}$ and $x$. The original argument~\cite{harperrmflowmoments} used $\sqrt{x}$, but this was changed to $x^{3/4}$ in the later high moments paper~\cite{harperrmfhigh} because that makes the next step of smoothing the outer sum to an integral rather cleaner. Although this effect is invisible at the level of detail in this survey, it may be of interest to a reader trying to master these arguments.} of \eqref{condhold}.

It is perhaps worth noting that Khintchine's inequality is not in fact a very deep statement here. If $a_p$ are any coefficients and $0 \leq q \leq 1$, then by H\"older's inequality we always have
$$ \E |\sum_{x^{3/4} < p \leq x} f(p) a_p|^{2q} \geq \frac{(\E |\sum_{x^{3/4} < p \leq x} f(p) a_p|^{2})^{2-q}}{(\E |\sum_{x^{3/4} < p \leq x} f(p) a_p|^{4})^{1-q}} = \frac{(\sum_{x^{3/4} < p \leq x} |a_p|^{2})^{2-q}}{(\E |\sum_{x^{3/4} < p \leq x} f(p) a_p|^{4})^{1-q}} . $$
Simply expanding the fourth power, we see $\E |\sum_{x^{3/4} < p \leq x} f(p) a_p|^{4}$ in the denominator is $= \sum_{x^{3/4} < p_1, ..., p_4 \leq x} a_{p_1} ... \overline{a_{p_4}} \textbf{1}_{p_1 p_2 = p_3 p_4} \ll \sum_{x^{3/4} < p_1, p_2 \leq x} |a_{p_1}|^2 |a_{p_2}|^2 = (\sum_{x^{3/4} < p \leq x} |a_p|^{2})^2$, giving the lower bound $\E |\sum_{x^{3/4} < p \leq x} f(p) a_p|^{2q} \gg (\sum_{x^{3/4} < p \leq x} |a_p|^{2})^q$.

\item When $q \geq 1$, H\"older's inequality no longer allows one to simply upper bound general quantities of the shape $\E X^{2q}$ by $(\E X^2)^q$. As a substitute, if we wish to upper bound $\E|\sum_{n \leq N} a_n f(n)|^{2q}$ for some coefficients $a_n \in \C$ (which in practice will be sums involving values of $f(m)$ on which we have conditioned), we can use a simple {\em hypercontractive inequality}. Thus for any real $q \geq 1$, we have
$$ \E|\sum_{n \leq N} a_n f(n)|^{2q} \leq \left( \sum_{n \leq N} |a_n|^2 d_{\lceil q \rceil}(n) \right)^q , $$
where $d_k(\cdot)$ denotes the $k$-fold divisor function (i.e. the number of $k$-tuples of natural numbers whose product is $\cdot$), and $\lceil q \rceil$ denotes the ceiling of $q$. Again, this inequality is not very deep, it is easily proved using H\"older's inequality to move to the case of integer $q$, and then expanding the $2q$-th power.

Because of the divisor function terms $d_{\lceil q \rceil}(n)$ in the hypercontractive inequality, one needs to be careful in constructing the sums to which it is applied, to avoid incurring losses. For example, applying this directly to $\E|\sum_{n \leq x} f(n)|^{2q}$ would yield an upper bound $\left( \sum_{n \leq x} d_{\lceil q \rceil}(n) \right)^q \asymp_q x^q \log^{q(\lceil q \rceil - 1)}x$, missing the truth in Theorem \ref{thmsthigh} by a factor $\asymp_q \log^{O(q)}x$. Likewise, applying it (after conditioning on $(f(p))_{p \leq P}$) to the first sum in \eqref{rsdisplay} would produce a contribution $\E \left( \sum_{\substack{P < m \leq x, \\ m \; \text{is} \; P \; \text{rough}}} d_{\lceil q \rceil}(m) |\sum_{\substack{n \leq x/m, \\ n \; \text{is} \; P \; \text{smooth}}} f(n)|^2 \right)^q$. Roughly speaking, without the coefficients $d_{\lceil q \rceil}(m)$ one expects the sum over $P$-rough numbers to give rise to a $1/\log P$ factor inside the bracket (as in our earlier description of the low moments case), but with these coefficients one expects a factor like $\frac{((\log x)/\log P)^{O(q)}}{\log P}$. Since we must certainly have $(\log x)/\log P \rightarrow \infty$ in order for the term $\sum_{\substack{n \leq x, \\ n \; \text{is} \; P \; \text{smooth}}} f(n)$ in \eqref{rsdisplay} to contribute negligibly, (and it turns out that when $q$ is large we actually need $P$ rather smaller than $x^{1/q}$), the extra $((\log x)/\log P)^{O(q)}$ would be fatal to obtaining sharp bounds.

To resolve this problem, one can replace \eqref{rsdisplay} by a more elaborate decomposition, e.g. writing
\begin{eqnarray}
\sum_{n \leq x} f(n) & = & \sum_{k=1}^{K} \sum_{\substack{n \leq x, \\ P_k < P(n) \leq P_{k-1}}} f(n) + \sum_{\substack{n \leq x, \\ n \; \text{is} \; P_K \; \text{smooth}}} f(n) \nonumber \\
& = & \sum_{k=1}^{K} \sum_{\substack{P_k < m \leq x, \\ p|m \; \Rightarrow \; P_k < p \leq P_{k-1}}} f(m) \sum_{\substack{n \leq x/m, \\ n \; \text{is} \; P_k \; \text{smooth}}} f(n) + \sum_{\substack{n \leq x, \\ n \; \text{is} \; P_K \; \text{smooth}}} f(n) , \nonumber 
\end{eqnarray}
for a suitable sequence $x = P_0 > P_1 > ... > P_K$. Recall the discussion in the footnote \footref{fnsubdiv}. If one takes $P_k = x^{e^{-k}}$, say, then when applying the hypercontractive inequality to a term $\sum_{\substack{P_k < m \leq x, \\ p|m \; \Rightarrow \; P_k < p \leq P_{k-1}}} f(m) \sum_{\substack{n \leq x/m, \\ n \; \text{is} \; P_k \; \text{smooth}}} f(n)$ one only incurs a loss of the shape $((\log P_{k-1})/\log P_k)^{O(q)} = e^{O(q)}$, which is acceptable for Theorem \ref{thmsthigh}.

With these modifications, when seeking moment upper bounds with $q \geq 1$ one can more or less follow the strategy from the low moments case, and end up needing to bound expressions like $x^q \E \left( \frac{1}{\log P_k} \int_{1}^{x/P_k} |\sum_{\substack{n \leq z, \\ n \; \text{is} \; P_k \; \text{smooth}}} f(n)|^{2} \frac{dz}{z^2} \right)^q$. One could immediately extend the integral to $\infty$ and apply Parseval's identity, as in the low moments case. But this would be inefficient for large $q$, because it would turn out that most of the contribution to the expectation came from the extraneous $z > x/P_k$. To mitigate this effect, one can use {\em Rankin's trick}, noting first that $\int_{1}^{x/P_k} |\sum_{\substack{n \leq z, \\ P_k \; \text{smooth}}} f(n)|^{2} \frac{dz}{z^2} \leq (\frac{x}{P_k})^{2q/\log x} \int_{1}^{x/P_k} |\sum_{\substack{n \leq z, \\ P_k \; \text{smooth}}} f(n)|^{2} \frac{dz}{z^{2+2q/\log x}} \leq e^{2q} \int_{1}^{x/P_k} |\sum_{\substack{n \leq z, \\ P_k \; \text{smooth}}} f(n)|^{2} \frac{dz}{z^{2+2q/\log x}}$. The $e^{2q}$ here contributes an acceptable factor $e^{O(q^2)}$ in Theorem \ref{thmsthigh}, and the shift by $2q/\log x$ helps to dampen the contribution from large $z$ if we now extend to an infinite integral.

In fact, Parseval's identity implies that
$$ \E \Biggl( \frac{1}{\log P_k} \int_{1}^{\infty} |\sum_{\substack{n \leq z, \\ n \; \text{is} \; P_k \; \text{smooth}}} f(n)|^{2} \frac{dz}{z^{2+2q/\log x}} \Biggr)^q =  \E \Biggl( \frac{1}{2\pi \log P_k} \int_{-\infty}^{\infty} \frac{|F_{P_k}(1/2 + \frac{q}{\log x} +it)|^{2}}{|1/2 + \frac{q}{\log x} + it|^2} dt \Biggr)^q , $$
so to prove (the upper bound part of) Theorem \ref{thmsthigh} it will suffice to bound $\E \left( \frac{1}{\log P} \int_{-1/2}^{1/2} |F_{P}(1/2 + \frac{q}{\log x} +it)|^{2} dt \right)^q$, on a suitable range of $P$ ($=P_k$).

Note that the Rankin shift by $2q/\log x$ manifested itself in the Euler product being evaluated at $1/2 + \frac{q}{\log x} +it$ rather than $1/2 + it$. As we will discuss further in section \ref{sechigh}, this shift roughly means that the contribution to the product from any primes $> e^{(\log x)/q} = x^{1/q}$ becomes (stochastically) bounded. This is very reasonable, since the $q$-th power of a prime larger than $x^{1/q}$ would exceed $x$, so we shouldn't expect it to contribute in a sharp bound for $\E|\sum_{n \leq x} f(n)|^{2q}$.

\item For lower bounds when $q \geq 1$, the argument from the low moments case ($0 \leq q \leq 1$) extends directly, because we still have $\E^{x^{3/4}} |\sum_{x^{3/4} < p \leq x} f(p) \sum_{n \leq x/p} f(n)|^{2q} \gg \left( \sum_{x^{3/4} < p \leq x} |\sum_{n \leq x/p} f(n)|^2 \right)^q$. In fact this follows simply from H\"older's inequality when $q \geq 1$, with no need for Khintchine's inequality. After the smoothing steps, one needs to introduce a similar Rankin shift as in the upper bound argument to discard the surplus from extending the $z$ integral to $\infty$, so again a quantity like $\E \left( \frac{1}{\log P} \int_{-1/2}^{1/2} |F_{P}(1/2 + \frac{q}{\log x} +it)|^{2} dt \right)^q$ emerges (with $P=x^{3/4}$, say).

\end{itemize}

\section{Low moments via barrier events}\label{seclow}
In this section, we shall discuss some of the further ideas involved in proving Theorem \ref{thmstlow}, beginning from the position we arrived at in section \ref{secep}. It would essentially suffice to show that uniformly for all large $P$ and all $2/3 \leq q \leq 1$, we have
\begin{equation}\label{chaostarget}
\E \Biggl( \int_{-1/2}^{1/2} |F_{P}(1/2+it)|^{2} dt \Biggr)^q \asymp \left( \frac{\log P}{1 + (1-q)\sqrt{\log\log P}} \right)^q .
\end{equation}
Recall that the relevant values of $P$ were $x^{1/\log\log x}$ (say) for the upper bound, and $x^{3/4}$ (say) for the lower bound, both of which satisfy $\log\log P \asymp \log\log x$.

\vspace{12pt}
It turns out that the left hand side of \eqref{chaostarget} is closely connected to a probabilistic object called {\em (critical) multiplicative chaos}. This is a fascinating and very active subject (along with the related topic of {\em log-correlated random processes}), motivating lots of our work, and is discussed in some detail in the introduction of \cite{harperrmflowmoments} along with further references. However, ultimately one can prove \eqref{chaostarget} in a fairly ``down to earth'' way, and this is how we shall try to present things here. In particular, we wish to motivate and explain the crucial role of so-called {\em barrier events}. For a short proof of the upper bound in \eqref{chaostarget} (for fixed $q$), but depending on rather a lot from the theory of multiplicative chaos, see Gorodetsky and Wong~\cite{gorwong}.

\vspace{12pt}
It seems most instructive to begin with the lower bound problem. A general strategy for lower bounding a quantity $\E I(f)^q$, where $I(f)$ is non-negative and $q \leq 1$, is to obtain a lower bound (or asymptotic) for $\E I(f)$ and an upper bound for some higher moment, say for $\E I(f)^2$. Then H\"older's inequality implies that
$$ \E I(f) = \E I(f)^{\frac{q}{2-q}} I(f)^{\frac{2(1-q)}{2-q}} \leq (\E I(f)^{q})^{\frac{1}{2-q}} (\E I(f)^2)^{\frac{1-q}{2-q}} , \;\;\;\;\; \text{and so} \;\;\;\;\; \E I(f)^q \geq \frac{(\E I(f))^{2-q}}{(\E I(f)^2)^{1-q}} . $$
(The reader may compare with our discussion of Khintchine's inequality, towards the end of section \ref{secep}.) Qualitatively, this simply reflects the fact that if the moments don't grow too rapidly as the power increases, then a significant portion of their size must come from fairly probable events (as opposed to rare, extreme events), and a good lower bound for low moments follows. For example, if we had $\E I(f) \gg C$ and $\E I(f)^2 \ll C^2$ (the strongest possible upper bound, in view of the Cauchy--Schwarz inequality), we would deduce a best possible lower bound $\E I(f)^q \gg C^q$.

If we try to apply this directly with $I(f) = \int_{-1/2}^{1/2} |F_{P}(1/2+it)|^{2} dt$, then the quantities we need to compute are
$$ \E I(f) = \int_{-1/2}^{1/2} \E |F_{P}(1/2+it)|^{2} dt , \;\;\; \text{and} \;\;\; \E I(f)^2 = \int_{-1/2}^{1/2} \int_{-1/2}^{1/2} \E |F_{P}(1/2+it)|^{2}|F_{P}(1/2+iu)|^{2} dt du . $$
Thanks to translation invariance in law, we can simplify these expressions by observing that $\E |F_{P}(1/2+it)|^{2} = \E |F_{P}(1/2)|^{2}$ for all $t \in \R$, and $\E |F_{P}(1/2+it)|^{2}|F_{P}(1/2+iu)|^{2} = \E |F_{P}(1/2)|^{2}|F_{P}(1/2+i(u-t))|^{2}$ for all $t,u \in \R$. Since $F_{P}(s)$ is an Euler product of independent factors, these expectations are not hard to calculate, and one finds that
\begin{equation}\label{onetwopoint}
\E |F_{P}(1/2)|^{2} \asymp \exp\{\sum_{p \leq P} \frac{1}{p}\} \asymp \log P ,
\end{equation}
$$ \E |F_{P}(1/2)|^{2}|F_{P}(1/2+ih)|^{2} \asymp \exp\{\sum_{p \leq P} \frac{2 + 2\cos(h\log p)}{p}\} \asymp \log^{2}P (\min\{\log P, \frac{1}{|h|}\})^2 , \;\;\; |h| \leq 1 . $$
Unfortunately, these estimates imply that $\E I(f) \asymp \log P$ but $\E I(f)^2 \asymp \log^{3}P$ (rather than $\log^{2}P$), giving a poor lower bound $\E I(f)^q \gg \frac{\log^{2-q}P}{\log^{3(1-q)}P} = \log^{2q-1}P$. Given the shape of Theorem \ref{thmstlow}, this failure should be unsurprising, since the bound we are actually seeking is not $\log^{q}P$ (as directly suggested by $\E I(f)$) but $\left( \frac{\log P}{1 + (1-q)\sqrt{\log\log P}} \right)^q$.

\vspace{12pt}
To improve the situation, we should think about possible inefficiencies in the above argument. The basic problem is that $\E I(f)^2$ is much larger than we might hope, because $\E |F_{P}(1/2)|^{2}|F_{P}(1/2+ih)|^{2}$ is too large when $|h|$ is small. This cannot be changed with the given obvious choice of $I(f)$, but since we ultimately just want a lower bound for $\E \left( \int_{-1/2}^{1/2} |F_{P}(1/2+it)|^{2} dt \right)^q$, {\em we are free to replace $I(f)$ from the beginning by any lower bound for $\int_{-1/2}^{1/2} |F_{P}(1/2+it)|^{2} dt$}.

A sensible choice of minorant for $\int_{-1/2}^{1/2} |F_{P}(1/2+it)|^{2} dt$ is not obvious, but we can get some idea about this by revisiting our estimate for $\E |F_{P}(1/2)|^{2}|F_{P}(1/2+ih)|^{2}$. The problematic factor $(\min\{\log P, \frac{1}{|h|}\})^2$ in that estimate arises from $\exp\{\sum_{p \leq P} \frac{2\cos(h\log p)}{p}\}$. The size may be explained by noting that when $p \leq e^{1/|h|}$, we have $\cos(h\log p) \approx 1$, whereas on larger $p$ we get oscillation and cancellation amongst the terms $\frac{2\cos(h\log p)}{p}$. Thus $\exp\{\sum_{p \leq P} \frac{2\cos(h\log p)}{p}\} \approx \exp\{\sum_{p \leq \min\{P,e^{1/|h|}\}} \frac{2}{p}\} \asymp (\min\{\log P, \frac{1}{|h|}\})^2$. This makes visible that depending on the size of $h$, different subproducts of $F_{P}(1/2), F_{P}(1/2+ih)$ are either highly correlated (producing blow-up) or relatively uncorrelated, namely the subproducts $F_{\min\{P,e^{1/|h|}\}}(1/2), F_{\min\{P,e^{1/|h|}\}}(1/2+ih)$ up to $\min\{P,e^{1/|h|}\}$ are highly correlated. We then arrive at a key idea in the area, that {\em rather than working with $\int_{-1/2}^{1/2} |F_{P}(1/2+it)|^{2} dt$ one should try to work with $\int_{-1/2}^{1/2} \textbf{1}_{B(t)} |F_{P}(1/2+it)|^{2} dt$}, where $B(t)$ is some {\em barrier event} that places restrictions on the sizes of {\em various subproducts} of $F_{P}(1/2+it)$.

To advance this idea, we must determine a sensible choice of $B(t)$. One can try to get a feel for this by thinking about what distributional behaviour of the random Euler products causes $\E |F_{P}(1/2)|^{2}$ to have size $\asymp \log P$, and causes $\E |F_{P}(1/2)|^{2}|F_{P}(1/2+ih)|^{2}$ to have size $\asymp \log^{2}P (\min\{\log P, \frac{1}{|h|}\})^2$. Beginning with the former, we have
$$ |F_{P}(1/2)|^{2} = \exp\{2\log|F_{P}(1/2)| \} = \exp\{2\Re \log F_{P}(1/2) \} = \exp\{-2 \sum_{p \leq P} \Re \log(1 - \frac{f(p)}{p^{1/2}}) \} . $$
The random summands $\Re \log(1 - \frac{f(p)}{p^{1/2}})$ here are independent (because the $f(p)$ are independent), and (keeping in mind the Taylor expansion $\log(1 - \frac{f(p)}{p^{1/2}}) = - \sum_{k=1}^{\infty} \frac{f(p)^k}{k p^{k/2}}$) it is easy to calculate that they have mean zero and variance $\E (\frac{\Re f(p)}{p^{1/2}})^2 + O(\frac{1}{p^{3/2}}) = \frac{1}{2p} + O(\frac{1}{p^{3/2}})$. So, in view of the central limit theorem for sums of independent random variables, we may expect $|F_{P}(1/2)|^{2}$ to have similar distributional behaviour as $\exp\{2 G_{P}\}$, where $G_P$ is a {\em Gaussian random variable} with mean zero and variance $\sum_{p \leq P} (\frac{1}{2p} + O(\frac{1}{p^{3/2}})) \sim (1/2)\log\log P$. Assuming for simplicity that the variance is exactly $(1/2)\log\log P$, an explicit calculation with the Gaussian probability density function shows that
\begin{eqnarray}
\E \exp\{2 G_{P}\} = \frac{1}{\sqrt{\pi \log\log P}} \int_{-\infty}^{\infty} e^{2z} e^{-z^{2}/\log\log P} dz & = & \frac{e^{\log\log P}}{\sqrt{\pi \log\log P}} \int_{-\infty}^{\infty} e^{-(z-\log\log P)^{2}/\log\log P} dz \nonumber \\
& = & \frac{e^{\log\log P}}{\sqrt{2\pi}} \int_{-\infty}^{\infty} e^{-w^{2}/2} dw = \log P . \nonumber
\end{eqnarray}
This matches with \eqref{onetwopoint}, but more importantly (looking at the values of $z$ that make the major contribution to the integrals) it suggests this size is produced by values of $G_{P}$ or $\log|F_{P}(1/2)|$ that are $\approx \log\log P$ (in fact in an interval of size $\asymp \sqrt{\log\log P}$ around $\log\log P$).

Proceeding similarly with $\E |F_{P}(1/2)|^{2}|F_{P}(1/2+ih)|^{2}$, we already noted that the blow-up there is created by the highly correlated subproducts $F_{\min\{P,e^{1/|h|}\}}(1/2), F_{\min\{P,e^{1/|h|}\}}(1/2+ih)$, so what we really need to analyse is $\E |F_{\min\{P,e^{1/|h|}\}}(1/2)|^2 |F_{\min\{P,e^{1/|h|}\}}(1/2+ih)|^2 \approx \E |F_{\min\{P,e^{1/|h|}\}}(1/2)|^4$. If we first look specifically at $\E |F_{P}(1/2)|^4$ (since this is easiest to write!), we may expect this to behave like
\begin{eqnarray}
\E \exp\{4 G_{P}\} = \frac{1}{\sqrt{\pi \log\log P}} \int_{-\infty}^{\infty} e^{4z} e^{-z^{2}/\log\log P} dz & = & \frac{e^{4\log\log P}}{\sqrt{\pi \log\log P}} \int_{-\infty}^{\infty} e^{-(z-2\log\log P)^{2}/\log\log P} dz \nonumber \\
& = & \frac{e^{4\log\log P}}{\sqrt{2\pi}} \int_{-\infty}^{\infty} e^{-w^{2}/2} dw = \log^{4}P . \nonumber
\end{eqnarray}
Again this agrees with our earlier observation that $\E |F_{P}(1/2)|^{4} \asymp \log^{4}P$, and we infer this size is produced by values of $\log|F_{P}(1/2)|$ that are $\approx 2\log\log P$. In particular, the size of $\log|F_{P}(1/2)|$ that produces the blow-up is {\em significantly larger} than the size that substantially contributes to $\E |F_{P}(1/2)|^{2}$. More generally, the size of $\log|F_{\min\{P,e^{1/|h|}\}}(1/2)|$ that substantially contributes to the blow-up term $(\min\{\log P, \frac{1}{|h|}\})^2$ in $\E |F_{P}(1/2)|^{2}|F_{P}(1/2+ih)|^{2}$ will be $\approx 2\log(\min\{\log P, \frac{1}{|h|}\})$, much larger than the size $\approx \log(\min\{\log P, \frac{1}{|h|}\})$ that substantially contributes to $\E |F_{\min\{P,e^{1/|h|}\}}(1/2)|^2$.

Motivated by all this, it seems reasonable to try working with $I(f) := \int_{-1/2}^{1/2} \textbf{1}_{B(t)} |F_{P}(1/2+it)|^{2} dt$, where $B(t)$ is the event that $\log|F_{y}(1/2+it)| \leq \log\log y$ (say) for all $2 \leq y \leq P$. We hope that this barrier will not reduce the size of $\E I(f)$ too much compared with the $\log P$ we had before, because such a condition still permits the sizes of subproducts that substantially contribute to $\E|F_{y}(1/2+it)|^2$; but that it will greatly reduce the size of $\E I(f)^2$, by forbidding the larger sizes of subproducts that would inflate this.

\vspace{12pt}
It turns out that, roughly speaking, such a strategy succeeds\footnote{It seems worth emphasising that the heuristics and motivation given up to this point certainly do not guarantee success, even assuming that any technical issues arising can be resolved. If we have the statements of Theorem \ref{thmstlow} and \eqref{chaostarget} in advance, we can be more confident, because we know that (a) the exponent in the moments should grow linearly with $q < 1$, suggesting the general strategy of comparing $\E I(f)$ and $\E I(f)^2$ might be appropriate and efficient; and (b) the size in \eqref{chaostarget} should not differ too much from the first guess $\log^{q}P$, so a choice of barrier $B(t)$ that doesn't alter $\E I(f)$ much should be our target. If we were considering the problem (or a related one) completely from scratch, more careful thought about the distribution of $F_{P}(1/2+it)$, and especially the interactions between the products for different $t$, might be needed.}, and we end this discussion with a few details of how the argument may actually be implemented. But first, we wish to flag up the ultimate source of the factor $\sqrt{\log\log P}$ in \eqref{chaostarget} (and thus the factor $\sqrt{\log\log x}$ in Theorem \ref{thmstlow}), which seems to be far the most interesting and a priori unexpected feature of the result.

The key point is that although imposing a restriction like $\log|F_{y}(1/2+it)| \leq \log\log y$ will not much alter the size of $\E|F_{y}(1/2+it)|^2$ (or of the full product $\E |F_{P}(1/2+it)|^{2}$) for any single $y$, when one imposes this restriction {\em simultaneously} for all (or many) $2 \leq y \leq P$ then it does reduce the expected size a bit. Indeed, setting $t=0$ to simplify the writing, the values $(\log|F_{y}(1/2)|)_{2 \leq y \leq P} = (-\sum_{p \leq y} \Re \log(1 - \frac{f(p)}{p^{1/2}}))_{2 \leq y \leq P}$ behave approximately like a {\em Gaussian random walk} in $y$, where the effective number of steps is the variance of the sum up to $P$, namely $\asymp \log\log P$. Then $\E |F_{P}(1/2)|^{2} \textbf{1}_{\log|F_{y}(1/2)| \leq \log\log y \; \forall \; 2 \leq y \leq P}$ is
$$ \asymp \log P \frac{\E |F_{P}(1/2)|^{2} \textbf{1}_{\log|F_{y}(1/2)| \leq \log\log y \; \forall \; 2 \leq y \leq P}}{\E |F_{P}(1/2)|^{2}} = \log P \frac{\E e^{2\log|F_{P}(1/2)|} \textbf{1}_{\log|F_{y}(1/2)| \leq \log\log y \; \forall \; 2 \leq y \leq P}}{\E e^{2\log|F_{P}(1/2)|}} . $$
The ratio is an {\em exponentially tilted probability}, i.e. the expectation of the event that $\log|F_{y}(1/2)| \leq \log\log y \; \forall \; 2 \leq y \leq P$ under the modified probability measure where everything is weighted by the random exponential factor $e^{2\log|F_{P}(1/2)|}$. If everything were exactly Gaussian, the very useful {\em Girsanov's theorem} (which really just involves completing the square in a computation with the Gaussian density) would imply that the tilted probability equals the probability that the same Gaussian random walk satisfies a certain modified condition. That condition would be roughly that the random walk remains below 0 (rather than $\log\log y$) for all $2 \leq y \leq P$. It turns out that the logarithms of our Euler products are sufficiently close to Gaussian that all this can be carried through for them as well. Finally, the classical probabilistic {\em Ballot Theorem} implies this probability is $\asymp \frac{1}{\sqrt{\log\log P}}$, with $\log\log P$ corresponding to the number of ``steps'' in the walk.

\vspace{12pt}
To actually prove the lower bound in \eqref{chaostarget}, it more or less suffices to follow the strategy just outlined, with $B(t)$ being the event that $\log|F_{y}(1/2+it)| \leq \log\log y + \min\{\sqrt{\log\log P}, \frac{1}{1-q}\} \; \forall \; 2 \leq y \leq P$. A Girsanov type calculation and the Ballot Theorem ultimately imply (note the extra term $\min\{\sqrt{\log\log P}, \frac{1}{1-q}\}$ in our barrier) that
$$ \E I(f) = \int_{-1/2}^{1/2} \E \textbf{1}_{B(t)} |F_{P}(1/2+it)|^{2} dt \asymp \log P \frac{\min\{\sqrt{\log\log P}, \frac{1}{1-q}\}}{\sqrt{\log\log P}} \asymp \frac{\log P}{1 + (1-q)\sqrt{\log\log P}} , $$
and we have (roughly speaking, see below)
\begin{eqnarray}
&& \E I(f)^2 = \int_{-1/2}^{1/2} \int_{-1/2}^{1/2} \E \textbf{1}_{B(t)} |F_{P}(1/2+it)|^{2} \textbf{1}_{B(u)} |F_{P}(1/2+iu)|^{2} dt du \nonumber \\
& \leq & \int_{-1}^{1} \E \textbf{1}_{B(0)} |F_{P}(1/2)|^{2} \textbf{1}_{B(h)} |F_{P}(1/2+ih)|^{2} dh \ll e^{2\min\{\sqrt{\log\log P}, \frac{1}{1-q}\}} \left(\frac{\log P}{1 + (1-q)\sqrt{\log\log P}} \right)^2 . \nonumber
\end{eqnarray}
Combining these bounds using H\"older's inequality, in the manner described earlier, gives
$$ \E \Biggl( \int_{-1/2}^{1/2} |F_{P}(1/2+it)|^{2} dt \Biggr)^q \geq \E I(f)^q \geq \frac{(\E I(f))^{2-q}}{(\E I(f)^2)^{1-q}} \gg \left(\frac{\log P}{1 + (1-q)\sqrt{\log\log P}} \right)^q . $$
Notice that the undesirable factor $e^{2\min\{\sqrt{\log\log P}, \frac{1}{1-q}\}}$, which is produced by the term $\min\{\sqrt{\log\log P}, \frac{1}{1-q}\}$ in our barrier, is killed off by the exponent $1-q$ to which we raise $\E I(f)^2$ (whereas a larger factor would not be). Most of the calculation required in the proof comes in showing that $\E \textbf{1}_{B(t)} |F_{P}(1/2+it)|^{2}$ and $\E \textbf{1}_{B(0)} |F_{P}(1/2)|^{2} \textbf{1}_{B(h)} |F_{P}(1/2+ih)|^{2}$ behave close to the way they would in the Gaussian case, which boils down to characteristic function calculations and careful use of the (two-dimensional) Berry--Esseen inequality.

There is one further technical but interesting point that deserves attention. When bounding $\E \textbf{1}_{B(0)} |F_{P}(1/2)|^{2} \textbf{1}_{B(h)} |F_{P}(1/2+ih)|^{2}$, one uses the barrier condition $B(h)$ to control the subproduct $F_{\min\{P,e^{1/|h|}\}}(1/2+ih)$, which as currently formulated would give an upper bound
$$ \leq e^{2\min\{\sqrt{\log\log P}, \frac{1}{1-q}\}} (\min\{\log P,1/|h|\})^2 \E \textbf{1}_{B(0)} |F_{P}(1/2)|^{2} \textbf{1}_{B(h)} |\prod_{\min\{P,e^{1/|h|}\} < p \leq P} (1 - \frac{f(p)}{p^{1/2+ih}})^{-1}|^{2} . $$
The products $F_{P}(1/2), \prod_{\min\{P,e^{1/|h|}\} < p \leq P} (1 - \frac{f(p)}{p^{1/2+ih}})^{-1}$ are then sufficiently uncorrelated that, roughly speaking\footnote{In fact, the parts of the barrier conditions $B(0), B(h)$ dealing with the subproducts up to $F_{\min\{P,e^{1/|h|}\}}(1/2), F_{\min\{P,e^{1/|h|}\}}(1/2+ih)$ remain highly correlated, but this can be understood and handled without too much difficulty (especially since we only need upper bounds at this point).}, the expectation on the previous line factors as
\begin{eqnarray}
& \approx & (\E \textbf{1}_{B(0)} |F_{P}(1/2)|^{2}) \cdot (\E \textbf{1}_{B(h)} |\prod_{\min\{P,e^{1/|h|}\} < p \leq P} (1 - \frac{f(p)}{p^{1/2+ih}})^{-1}|^{2}) \nonumber \\
& \approx & \left(\frac{\log P}{1 + (1-q)\sqrt{\log\log P}} \right) \cdot \left(\frac{\max\{|h|\log P, 1\}}{1 + (1-q)\sqrt{\log\log P}} \right) . \nonumber
\end{eqnarray}
The term $\max\{|h|\log P, 1\}$ here is $\E \prod_{\min\{P,e^{1/|h|}\} < p \leq P} |1 - \frac{f(p)}{p^{1/2+ih}}|^{-2}$. Multiplying by $e^{2\min\{\sqrt{\log\log P}, \frac{1}{1-q}\}} (\min\{\log P,1/|h|\})^2$ and integrating over $h$ would then deliver a bound
$$ \E I(f)^2 \lesssim e^{2\min\{\sqrt{\log\log P}, \frac{1}{1-q}\}} \left(\frac{\log P}{1 + (1-q)\sqrt{\log\log P}} \right)^2 \int_{-1}^{1} \min\{\log P,1/|h|\} dh . $$
Unfortunately, the integral over $h$ would still contribute an unwanted factor $\asymp \log\log P$ here, and so we would not quite obtain a satisfactory bound for $\E I(f)^2$.

To fix this, we actually make a slightly different choice of barrier, taking $B(t)$ to be the event that $\log|F_{y}(1/2+it)| \leq \log\log y - 2\log\log\log y + \min\{\sqrt{\log\log P}, \frac{1}{1-q}\} \; \forall \; 100 \leq y \leq P$ (say). Tracing everything through, the extra subtracted term $-2\log\log\log y$ means that we end up with $\int_{-1}^{1} \min\{\frac{\log P}{(\log\log P)^4},\frac{1}{|h| \log^{4}(2/|h|)}\} dh$ (which is bounded) rather than $\int_{-1}^{1} \min\{\log P,1/|h|\} dh$. And happily this modified barrier {\em does not change} the order of magnitude of $\E I(f)$, in fact there is lots of useful flexibility in this regard, we could adjust the barrier up or down by roughly $\sqrt{\log\log y}$ (the typical fluctuations of the random walk up to $y$) without altering $\E I(f)$.

\vspace{12pt}
To prove the upper bound in \eqref{chaostarget}, one can use barrier events of a similar shape as in the lower bound discussion. We may assume that $q \leq 1 - \frac{1}{\sqrt{\log\log P}}$, otherwise the desired upper bound is trivial. For a large parameter $C$, let $B_{C}(t)$ be the event that $\log|F_{y}(1/2+it)| \leq \log\log y + 2\log\log\log y + \frac{C}{1-q} \; \forall \; 100 \leq y \leq P$. Note the added term $2\log\log\log y$, as opposed to the subtracted $-2\log\log\log y$ in the lower bound argument. Also let $\mathcal{G}_{C}$ denote the ``good'' event that $B_{C}(t)$ holds for all $|t| \leq 1/2$. Then clearly $\E \left( \int_{-1/2}^{1/2} |F_{P}(1/2+it)|^{2} dt \right)^q$ is
\begin{eqnarray}
& = & \E \textbf{1}_{\mathcal{G}_{C} \; \text{holds}} \Biggl( \int_{-1/2}^{1/2} |F_{P}(1/2+it)|^{2} dt \Biggr)^q + \E \textbf{1}_{\mathcal{G}_{C} \; \text{fails}} \Biggl( \int_{-1/2}^{1/2} |F_{P}(1/2+it)|^{2} dt \Biggr)^q \nonumber \\
& \leq & \E \Biggl( \int_{-1/2}^{1/2} \textbf{1}_{B_{C}(t)} |F_{P}(1/2+it)|^{2} dt \Biggr)^q + \E \textbf{1}_{\mathcal{G}_{C} \; \text{fails}} \Biggl( \int_{-1/2}^{1/2} |F_{P}(1/2+it)|^{2} dt \Biggr)^q . \nonumber
\end{eqnarray}
By H\"older's inequality, the first integral is $\leq \left( \int_{-1/2}^{1/2} \E \textbf{1}_{B_{C}(t)} |F_{P}(1/2+it)|^{2} dt \right)^q$, and this is $\ll \left(\frac{C \log P}{(1-q)\sqrt{\log\log P}} \right)^q \leq C \left(\frac{\log P}{(1-q)\sqrt{\log\log P}} \right)^q$ by a Girsanov--Ballot Theorem calculation. Notice that, as discussed earlier, the added $2\log\log\log y$ in the definition of $B_{C}(t)$ makes no visible difference to the Girsanov--Ballot Theorem bound.

Unlike with lower bounds, we of course cannot just discard the second integral $\E \textbf{1}_{\mathcal{G}_{C} \; \text{fails}} \left( \int_{-1/2}^{1/2} |F_{P}(1/2+it)|^{2} dt \right)^q$. We explain how to handle this, in the style of a nice recent paper\footnote{The original argument of Harper~\cite{harperrmflowmoments} uses a sequence of applications of H\"older's inequality to replace the exponent $q$ by $q', q'', ...$, halving the distance to 1 at each step. This amounts to considering barriers with a sequence of different $C$ values (replacing $\frac{C}{1-q}$ by $\frac{C}{1-q'}, \frac{C}{1-q''}, ...$ is equivalent to replacing it by $\frac{C'}{1-q}, \frac{C''}{1-q}, ...$, for suitable $C', C'', ...$), as Soundararajan and Zaman~\cite{soundzam} do, but their direct presentation seems the simpler and clearer approach.} of Soundararajan and Zaman~\cite{soundzam}. If $C' > C$ is a further parameter, the integral is
\begin{eqnarray}
& = & \E \textbf{1}_{\mathcal{G}_{C} \; \text{fails but} \; \mathcal{G}_{C'} \; \text{holds}} \Biggl( \int_{-1/2}^{1/2} |F_{P}(1/2+it)|^{2} dt \Biggr)^q + \E \textbf{1}_{\mathcal{G}_{C'} \; \text{fails}} \Biggl( \int_{-1/2}^{1/2} |F_{P}(1/2+it)|^{2} dt \Biggr)^q \nonumber \\
& \leq & \E \textbf{1}_{\mathcal{G}_{C} \; \text{fails}} \Biggl( \int_{-1/2}^{1/2} \textbf{1}_{B_{C'}(t)} |F_{P}(1/2+it)|^{2} dt \Biggr)^q + \E \textbf{1}_{\mathcal{G}_{C'} \; \text{fails}} \Biggl( \int_{-1/2}^{1/2} |F_{P}(1/2+it)|^{2} dt \Biggr)^q . \nonumber
\end{eqnarray}
By H\"older's inequality, but now treating the factor $\textbf{1}_{\mathcal{G}_{C} \; \text{fails}}$ non-trivially, the first term here is $\leq \p(\mathcal{G}_{C} \; \text{fails})^{1-q} \cdot \left( \int_{-1/2}^{1/2} \E \textbf{1}_{B_{C'}(t)} |F_{P}(1/2+it)|^{2} dt \right)^q \ll \p(\mathcal{G}_{C} \; \text{fails})^{1-q} \cdot C' \left(\frac{\log P}{(1-q)\sqrt{\log\log P}} \right)^q$. To estimate $\p(\mathcal{G}_{C} \; \text{fails})$, note that for any given $100 \leq y \leq P$ and $|t| \leq 1/2$, we have
\begin{eqnarray}
&& \p(\log|F_{y}(1/2+it)| > \log\log y + 2\log\log\log y + \frac{C}{1-q}) \nonumber \\
& \leq & \frac{\E|F_{y}(1/2+it)|^2}{e^{2(\log\log y + 2\log\log\log y + \frac{C}{1-q})}} \ll e^{-\frac{2C}{1-q}} \frac{1}{\log y (\log\log y)^4} . \nonumber
\end{eqnarray}
Since the most rapidly oscillating terms $p^{-it} = e^{-it\log p}$ involved in $\log|F_{y}(1/2+it)|$ rotate with speed $\log p \leq \log y$, it turns out that one can control $F_{y}(1/2+it)$ for all $|t| \leq 1/2$ by controlling it at a {\em net of points} $t$ with slightly tighter spacing than $1/\log y$. For example, a net of $\ll \log y (\log\log y)$ points is sufficient. Furthermore, it suffices to handle values of $y$ of the shape $e^{e^j}$, say, so that $\log\log y$ increments by 1. Thus, roughly speaking, the union bound implies that
$$ \p(\mathcal{G}_{C} \; \text{fails}) \lesssim \sum_{y = e^{e^j} \leq P} \log y (\log\log y) e^{-\frac{2C}{1-q}} \frac{1}{\log y (\log\log y)^4} = e^{-\frac{2C}{1-q}} \sum_{y = e^{e^j} \leq P} \frac{1}{j^3} \ll e^{-\frac{2C}{1-q}} . $$
Notice how the added term $2\log\log\log y$ in the barrier ultimately led to this sum over $y$ being uniformly bounded. We now see that $\E \textbf{1}_{\mathcal{G}_{C} \; \text{fails}} \left( \int_{-1/2}^{1/2} |F_{P}(1/2+it)|^{2} dt \right)^q $ is
$$ \lesssim e^{-2C} C' \left(\frac{\log P}{(1-q)\sqrt{\log\log P}} \right)^q + \E \textbf{1}_{\mathcal{G}_{C'} \; \text{fails}} \Biggl( \int_{-1/2}^{1/2} |F_{P}(1/2+it)|^{2} dt \Biggr)^q . $$
Applying this argument repeatedly, with a sequence $C < C' < C'' < ...$ such that the sum of the terms $C, e^{-2C} C', e^{-2C'} C'', ...$ is uniformly bounded (e.g. the sequence of natural numbers would suffice), finishes the proof. Actually one can stop the argument as soon as the $C$ value exceeds $(1-q)\sqrt{\log\log P}$, since then the trivial bound $\E \left( \int_{-1/2}^{1/2} |F_{P}(1/2+it)|^{2} dt \right)^q \ll \log^{q}P$ is as good as the bound one hopes for with the barrier $B_{C}(t)$ present. 

As a final technical remark, we note that there are different possibilities for making the above sketch argument (i.e. the estimation of $\p(\mathcal{G}_{C} \; \text{fails})$) fully rigorous. The original paper of Harper~\cite{harperrmflowmoments} used a modified definition of $\mathcal{G}_{C}$, where from the start the barrier conditions were only required to hold at a net of points $t$. Then the calculation of $\p(\mathcal{G}_{C} \; \text{fails})$ can be performed exactly as described, but one works a little more in the Girsanov--Ballot Theorem calculations to see that only having the barrier at a point near $t$ still suffices to produce the Ballot Theorem saving. Soundararajan and Zaman~\cite{soundzam} do not modify $\mathcal{G}_{C}$, but then they must incorporate a further discretisation (Sobolev--Gallagher type) argument into their estimation of $\p(\mathcal{G}_{C} \; \text{fails})$.

\section{High moments via Euler product correlations}\label{sechigh}
To deduce Theorem \ref{thmsthigh}, beginning from the position reached in section \ref{secep}, it would essentially suffice to prove an estimate like
\begin{equation}\label{hightarget}
\E \Biggl( \int_{-1/2}^{1/2} |F_{P}(1/2+\frac{q}{\log x} + it)|^{2} dt \Biggr)^q = e^{O(q^2)} \left(\frac{\min\{\log P, (\log x)/q\}}{\log 2q}\right)^{q^2 - q + 1} ,
\end{equation}
uniformly on a suitable range of $P \leq x$. The required range of $P$ depends on $x$ and on $q$ (which we recall may be a growing function of $x$ in Theorem \ref{thmsthigh}), so for simplicity we give no details about it here, except in a couple of places where it becomes relevant to qualitative features of the overall bounds. Notice the factor $\frac{1}{q\log 2q}$ visible in the bracket on the right (when $P \geq x^{1/q}$), which is responsible for the term $e^{-q^{2}\log q - q^{2}\log\log(2q)}$ in the theorem.

The shift $q/\log x$ in the Euler product in \eqref{hightarget} means that the contribution from any primes $> x^{1/q}$ becomes stochastically bounded, in other words $|F_{P}(1/2+\frac{q}{\log x} + it)|$ usually behaves in roughly the same way as $|F_{\min\{P,x^{1/q}\}}(1/2+ it)|$. For example, using the independence of the $f(p)$, it is again fairly easy to calculate that
$$ \E|\prod_{\min\{P,x^{1/q}\} < p \leq P} (1 - \frac{f(p)}{p^{1/2+q/\log x+it}})^{-1}|^{2q} = \exp\{\sum_{\min\{P,x^{1/q}\} < p \leq P} \frac{q^2}{p^{1+2q/\log x}} + O(\frac{q^2}{\log(2q)}) \} . $$
(Strictly speaking, this is true provided that $\min\{P,x^{1/q}\} \geq 100q^2$, say.) If $P > x^{1/q}$ then we have $\sum_{x^{1/q} < p \leq P} \frac{q^2}{p^{1+2q/\log x}} \leq \frac{q^3}{\log x} \sum_{x^{1/q} < p \leq P} \frac{\log p}{p^{1+2q/\log x}} \ll \frac{q^3}{\log x} \int_{x^{1/q}}^{\infty} \frac{1}{w^{1+2q/\log x}} dw$, using e.g. the classical Chebychev estimates from prime number theory. Performing the integral, we see this is all $\ll q^2$, which would contribute an acceptable $e^{O(q^2)}$ to \eqref{hightarget}. So (replacing $P$ by $\min\{P,x^{1/q}\}$ in the Euler product, and then relabelling this as $P$ for simplicity) we may ignore the shift by $q/\log x$, and work as though $\min\{\log P, (\log x)/q\}$ is replaced by $\log P$ in our target bound on the right hand side of \eqref{hightarget}, provided we can do everything with sufficient uniformity in $P$. Notice that it is an effect of the ``large primes'' that causes $\min\{\log P, (\log x)/q\}$ to arise here, and ultimately contributes $\frac{1}{q}$ to the crucial factor $\frac{1}{q\log 2q}$ that we observed above.

\vspace{12pt}
A key observation is that as $q$ increases, the values of $|F_{P}(1/2+it)|$ that significantly contribute to $\E ( \int_{-1/2}^{1/2} |F_{P}(1/2+ it)|^{2} dt )^q$ are more extreme, larger, less probable values. This is of course a very general point, that we already made in the Introduction when discussing our overall bounds for $\E |\sum_{n \leq x} f(n)|^{2q}$, and many times throughout the discussion in section \ref{seclow}. More specifically, in our barrier constructions we noted that the values of $|F_{P}(1/2+it)|$ that make the major contribution to $\E \int_{-1/2}^{1/2} |F_{P}(1/2+ it)|^{2} dt \asymp \log P$ are those where $\log|F_{P}(1/2+it)| \approx \log\log P$. Such values are {\em just} rare enough that they are unlikely to actually occur for $|t| \leq 1/2$, hence the small reduction in the expected value when one inserts the barriers $B(t)$ and $B_{C}(t)$ (which are obeyed with high probability). This is all ultimately responsible for the subtle size mismatch between the $2q$-th moments for $q < 1$ and $q=1$. For $q > 1$, the important values of $|F_{P}(1/2+ it)|$ will then be even larger and more improbable, and we should expect {\em only a small number of random $t$-values} (and, by continuity, short intervals around them) to significantly contribute to $\E ( \int_{-1/2}^{1/2} |F_{P}(1/2+ it)|^{2} dt )^q$. This can guide the assembly of our arguments.

\vspace{12pt}
Armed with the above observations, proving a good lower bound for $\E ( \int_{-1/2}^{1/2} |F_{P}(1/2+ it)|^{2} dt )^q$ becomes quite straightforward. Firstly it is very convenient, to streamline the manipulation of fractional powers, to replace $\int_{-1/2}^{1/2}$ with a discrete sum. We can neatly achieve this using Jensen's inequality. Thus
\begin{eqnarray}\label{shortintdisp}
&& \E ( \int_{-1/2}^{1/2} |F_{P}(1/2+ it)|^{2} dt )^q = \frac{1}{\log^{q}P} \E ( \log P \int_{-1/2}^{1/2} |F_{P}(1/2+ it)|^{2} dt )^q \nonumber \\
& \approx & \frac{1}{\log^{q}P} \E \Biggl( \sum_{|k| \leq (\log P)/2} \log P \int_{-1/(2\log P)}^{1/(2\log P)} |F_{P}(1/2+ \frac{ik}{\log P} + it)|^{2} dt \Biggr)^q ,
\end{eqnarray}
and since $|F_{P}(1/2+ \frac{ik}{\log P} + it)|^{2} = e^{2\log|F_{P}(1/2+ \frac{ik}{\log P} + it)|}$ and the exponential function is convex, Jensen's inequality (applied to the normalised integral $\log P \int_{-1/(2\log P)}^{1/(2\log P)} dt$) implies this is all
$$ \geq \frac{1}{\log^{q}P} \E \Biggl( \sum_{|k| \leq (\log P)/2} e^{2\log P \int_{-1/(2\log P)}^{1/(2\log P)} \log|F_{P}(1/2+ \frac{ik}{\log P} + it)| dt} \Biggr)^q . $$
Here $e^{2\log P \int_{-1/(2\log P)}^{1/(2\log P)} \log|F_{P}(1/2+ \frac{ik}{\log P} + it)| dt}$ behaves in essentially the same way as $e^{2 \log|F_{P}(1/2+ \frac{ik}{\log P})|}$. For simplicity we shall write the rest of the argument for $\frac{1}{\log^{q}P} \E ( \sum_{|k| \leq (\log P)/2} e^{2 \log|F_{P}(1/2+ \frac{ik}{\log P})|} )^q = \frac{1}{\log^{q}P} \E ( \sum_{|k| \leq (\log P)/2} |F_{P}(1/2+ \frac{ik}{\log P})|^2 )^q$, but one can perform all the same calculations rigorously for $e^{2\log P \int_{-1/(2\log P)}^{1/(2\log P)} \log|F_{P}(1/2+ \frac{ik}{\log P} + it)| dt}$.

Since $q \geq 1$, we have
$$ \frac{1}{\log^{q}P} \E \Biggl( \sum_{|k| \leq (\log P)/2} |F_{P}(1/2+ \frac{ik}{\log P})|^2 \Biggr)^q \geq \frac{1}{\log^{q}P} \E \sum_{|k| \leq (\log P)/2} |F_{P}(1/2+ \frac{ik}{\log P})|^{2q} . $$
This step would be very wasteful if many of the products $|F_{P}(1/2+ \frac{ik}{\log P})|^2$ made substantial contributions to the sum, but we observed earlier that here we expect the dominant contribution to come from {\em just a few} large products (at some random $k$).

Finally, it only remains to estimate $\E |F_{P}(1/2+ \frac{ik}{\log P})|^{2q}$. Since $F_{P}(s)$ is an Euler product of independent factors, this expectation is not hard to calculate, and {\em provided that $P \geq 100q^2$ (say)} one finds
$$ \E |F_{P}(1/2+ \frac{ik}{\log P})|^{2q} = \exp\{\sum_{100q^2 < p \leq P} \frac{q^2}{p} + O(\frac{q^2}{\log(2q)}) \} = \Biggl(\frac{\log P}{\log(100q^2)} \Biggr)^{q^2} e^{O(\frac{q^2}{\log(2q)})} . $$
Notice that the primes $\leq 100q^2$ are handled separately here: their contribution to the expectation of the Euler product is bounded trivially and goes into the $O(\frac{q^2}{\log(2q)})$ term, they do {\em not} produce a larger contribution $\sum_{p \leq 100q^2} \frac{q^2}{p}$. This ``small primes'' effect (the breakdown of Gaussian tail behaviour on the small primes when looking at very high moments) is thus responsible for the denominator $\log(100q^2) \asymp \log 2q$ in \eqref{hightarget}.

Putting everything together, one has a lower bound $\frac{1}{\log^{q}P} \sum_{|k| \leq (\log P)/2} (\frac{\log P}{\log(100q^2)} )^{q^2} e^{O(\frac{q^2}{\log(2q)})} =  e^{O(q^2)} \left(\frac{\log P}{\log 2q}\right)^{q^2 - q + 1}$, as desired. We also remark that one source of the upper bound restriction $q \leq \frac{c\log x}{\log\log x}$ in Theorem \ref{thmsthigh} is the need to have something like $P \geq 100q^2$, where $P$ may have size around $x^{1/q}$.

\vspace{12pt}
Moving to upper bounds for $\E ( \int_{-1/2}^{1/2} |F_{P}(1/2+ it)|^{2} dt )^q$, again it is convenient to replace $\int_{-1/2}^{1/2}$ with a discrete sum, but now Jensen's inequality goes in the wrong direction. Instead, a simple application of H\"older's inequality to the normalised integral $\log P \int_{-1/(2\log P)}^{1/(2\log P)} dt$ implies, since $q \geq 1$, that \eqref{shortintdisp} is
$$ \leq \frac{1}{\log^{q}P} \log P \int_{-1/(2\log P)}^{1/(2\log P)} \E \Biggl( \sum_{|k| \leq (\log P)/2} |F_{P}(1/2+ \frac{ik}{\log P} + it)|^{2} \Biggr)^q dt . $$
Since the joint distribution of $(F_{P}(1/2+\frac{ik}{\log P}+it))_{|k| \leq (\log P)/2}$ is the same for all $t \in \R$, the expectation here is the same for all $t \in \R$, so the above expression is in fact $= \frac{1}{\log^{q}P} \E ( \sum_{|k| \leq (\log P)/2} |F_{P}(1/2+ \frac{ik}{\log P})|^{2} )^q$.

We can gain further insight by rewriting $\frac{1}{\log^{q}P} \E ( \sum_{|k| \leq (\log P)/2} |F_{P}(1/2+ \frac{ik}{\log P})|^{2} )^q$ as
\begin{equation}\label{breakqdisplay}
\frac{1}{\log^{q}P} \sum_{|k| \leq (\log P)/2} \E |F_{P}(1/2+ \frac{ik}{\log P})|^{2} \Biggl( \sum_{|m| \leq (\log P)/2} |F_{P}(1/2+ \frac{im}{\log P})|^{2} \Biggr)^{q-1} .
\end{equation}
Recall once more that since $q \geq 1$, we expect the dominant contribution to come from just a few large Euler products at some random $k$, reinforced by {\em the same products} inside the $(q-1)$-st power. We can make this quite vivid by computing the ``correlation'' of $|F_{P}(1/2+ \frac{ik}{\log P})|^{2}$ and $|F_{P}(1/2+ \frac{im}{\log P})|^{2(q-1)}$, which again (since these are products of independent factors) is fundamentally a straightforward computation. One finds, provided $P \geq 100q^2$, that
\begin{eqnarray}\label{highcorrcalc}
&& \E |F_{P}(1/2+ \frac{ik}{\log P})|^{2} |F_{P}(1/2+ \frac{im}{\log P})|^{2(q-1)} \nonumber \\
& = & \exp\Biggl\{\sum_{100q^2 < p \leq P} \frac{(1+(q-1)^2 + 2(q-1)\cos(\frac{(m-k)\log p}{\log P}))}{p} + O(\frac{q^2}{\log(2q)}) \Biggr\} \nonumber \\
& = & \Biggl(\frac{\log P}{\log(100q^2)} \Biggr)^{1 + (q-1)^2} \Biggl(\frac{\log P}{1+|m-k|} \Biggr)^{2(q-1)} e^{O(\frac{q^2}{\log(2q)})} .
\end{eqnarray}
When $m=k$ this has the size $(\frac{\log P}{\log(100q^2)} )^{q^2} e^{O(\frac{q^2}{\log(2q)})}$ that we observed in our discussion of lower bounds, but as $|m-k|$ increases the size goes down (increasingly rapidly as $q$ becomes larger).

When $q-1 \geq 1$, we have the option of applying H\"older's inequality to the $(q-1)$-st power in \eqref{breakqdisplay}. An immediate application, bounding this by $\ll (\log P)^{q-2} \sum_{|m| \leq (\log P)/2} |F_{P}(1/2+ \frac{im}{\log P})|^{2(q-1)}$, is inefficient--- recall that we expect only a bounded number of $m$ values near to $k$ to typically contribute, and the factor $(\log P)^{q-2}$ multiplying everything (including the $m=k$ term) does not reflect this. But with only slightly more ingenuity, we can succeed. For example, we may note that $( \sum_{|m| \leq (\log P)/2} |F_{P}(1/2+ \frac{im}{\log P})|^{2} )^{q-1}$ is
\begin{eqnarray}
& = & \Biggl( \sum_{|m| \leq (\log P)/2} \frac{(1+|m-k|)^2 |F_{P}(1/2+ \frac{im}{\log P})|^{2}}{(1+|m-k|)^2} \Biggr)^{q-1} \nonumber \\
& \leq & e^{O(q)} \sum_{|m| \leq (\log P)/2} \frac{1}{(1+|m-k|)^2} (1+|m-k|)^{2(q-1)} |F_{P}(1/2+ \frac{im}{\log P})|^{2(q-1)} , \nonumber
\end{eqnarray}
by applying H\"older's inequality to the counting measure weighted by $\frac{1}{(1+|m-k|)^2}$. If we then multiply by $|F_{P}(1/2+ \frac{ik}{\log P})|^{2}$ and take expectations, the decaying factor $(\frac{1}{1+|m-k|})^{2(q-1)}$ in \eqref{highcorrcalc} nullifies the factor $(1+|m-k|)^{2(q-1)}$ from our weighted application of H\"older's inequality, and we deduce $\E |F_{P}(1/2+ \frac{ik}{\log P})|^{2} ( \sum_{|m| \leq (\log P)/2} |F_{P}(1/2+ \frac{im}{\log P})|^{2} )^{q-1} \leq e^{O(\frac{q^2}{\log(2q)})} \sum_{|m| \leq (\log P)/2} \frac{1}{(1+|m-k|)^2} (\frac{\log P}{\log(100q^2)} )^{1 + (q-1)^2} (\log P)^{2(q-1)} = e^{O(\frac{q^2}{\log(2q)})} (\frac{\log P}{\log(100q^2)} )^{q^2}$. This implies that \eqref{breakqdisplay} is $\leq \frac{1}{\log^{q}P} \sum_{|k| \leq (\log P)/2} e^{O(\frac{q^2}{\log(2q)})} (\frac{\log P}{\log(100q^2)} )^{q^2} = e^{O(\frac{q^2}{\log(2q)})} (\frac{\log P}{\log(100q^2)} )^{q^2 - q +1}$, a sharp bound.

\vspace{12pt}
It only remains to prove a good upper bound when $1 < q < 2$, which in fact is the most challenging part of Theorem \ref{thmsthigh} (the hardest case of all being when $q=q(x)$ tends down to 1). The argument is quite technical to execute properly and we shall not present many details, see section 5.4 of the original paper~\cite{harperrmfhigh} for the full proof. Instead, we briefly describe the ideas and tools required to adapt the above (fairly simple) $q \geq 2$ argument to the range $1 < q < 2$.

Since $q-1 < 1$, we can no longer apply H\"older's inequality only to the sum over $m$ in \eqref{breakqdisplay} to deliver Euler products $|F_{P}(1/2+ \frac{im}{\log P})|^{2(q-1)}$ of the shape we expect. (Recall that we expect that whole sum, raised to the power $q-1$, to typically behave in roughly the same way as $|F_{P}(1/2+ \frac{ik}{\log P})|^{2(q-1)}$.) Instead, we look to craft a suitable application of H\"older's inequality to the full expectation $\E |F_{P}(1/2+ \frac{ik}{\log P})|^{2} ( \sum_{|m| \leq (\log P)/2} |F_{P}(1/2+ \frac{im}{\log P})|^{2} )^{q-1}$. The obvious approach is to raise $( \sum_{|m| \leq (\log P)/2} |F_{P}(1/2+ \frac{im}{\log P})|^{2} )^{q-1}$ to the power $1/(q-1) > 1$, so that the sum over $m$ is no longer trapped inside a fractional power. This would leave us raising $|F_{P}(1/2+ \frac{ik}{\log P})|^{2}$ to the complementary exponent $1/(2-q)$. But there are two clear reasons why such an argument cannot be efficient. Firstly, this completely decouples the sum over $m$ from the point $k$, whereas we expect the sum to be dominated by terms around $k$ precisely because of the multiplying product $|F_{P}(1/2+ \frac{ik}{\log P})|^{2}$. Secondly, if we expect a final answer roughly the same size as $\E |F_{P}(1/2+ \frac{ik}{\log P})|^{2} |F_{P}(1/2+ \frac{ik}{\log P})|^{2(q-1)} = \E |F_{P}(1/2+ \frac{ik}{\log P})|^{2q} \asymp \log^{q^2}P$, then (because the dependence on $q$ is not linear) we need to keep terms of roughly the shape $\E |F_{P}(1/2+ \frac{ik}{\log P})|^{2q}$ in both factors that emerge from H\"older's inequality. We see that $\E |F_{P}(1/2+ \frac{ik}{\log P})|^{2/(2-q)}$ is not of this shape.

We can progress by exploiting the product structure of $F_{P}(s)$, in a way inspired by the consideration of subproducts in the barrier arguments of section \ref{seclow}. Thus if we set $C = C(q) := e^{1/(q-1)}$, in the part of the sum over $m$ where $C^{d-1} \leq |m-k| \leq C^d$ (for some $d$) we expect the subproducts $F_{P^{1/C^d}}(1/2+ \frac{im}{\log P})$ to be highly correlated with $F_{P^{1/C^d}}(1/2+ \frac{ik}{\log P})$, whilst the subproducts over primes $P^{1/C^d} < p \leq P$ should be fairly uncorrelated. It then makes sense to separate $|F_{P}(1/2+ \frac{ik}{\log P})|^2$ as $|F_{P}(1/2+ \frac{ik}{\log P})|^2 = |F_{P^{1/C^d}}(1/2+ \frac{ik}{\log P})|^2 H_{P,d}(1/2+ \frac{ik}{\log P})$, so that we can potentially apply H\"older's inequality with different exponents attached to $|F_{P^{1/C^d}}(1/2+ \frac{ik}{\log P})|^2$ and to $H_{P,d}(1/2+ \frac{ik}{\log P}) := \prod_{P^{1/C^d} < p \leq P} |1 - \frac{f(p)}{p^{1/2+ik/\log P}}|^{-2}$.

As an initial attempt, and assuming that $k=0$ to simplify the writing, we may rewrite $\E |F_{P}(1/2)|^{2} ( \sum_{|m| \leq (\log P)/2} |F_{P}(1/2+ \frac{im}{\log P})|^{2} )^{q-1}$ as
\begin{eqnarray}\label{ddecomp}
& \approx & \E |F_{P}(1/2)|^{2} \Biggl( \sum_{d \leq (q-1)\log P + 1} \sum_{\substack{C^{d-1} \leq |m| \leq C^d, \\ |m| \leq (\log P)/2}} |F_{P}(1/2+ \frac{im}{\log P})|^{2} \Biggr)^{q-1} \nonumber \\
& \leq & \sum_{d \leq (q-1)\log P + 1}  \E |F_{P}(1/2)|^{2} \Biggl( \sum_{\substack{C^{d-1} \leq |m| \leq C^d, \\ |m| \leq (\log P)/2}} |F_{P}(1/2+ \frac{im}{\log P})|^{2} \Biggr)^{q-1} \\
& = & \sum_{d \leq (q-1)\log P + 1}  \E |F_{P^{1/C^d}}(\frac{1}{2})|^{2q(2-q)} H_{P,d}(\frac{1}{2}) \Biggl( \sum_{\substack{C^{d-1} \leq |m| \leq C^d, \\ |m| \leq (\log P)/2}} |F_{P^{1/C^d}}(\frac{1}{2})|^{2(q-1)} |F_{P}(\frac{1}{2}+ \frac{im}{\log P})|^{2} \Biggr)^{q-1} . \nonumber 
\end{eqnarray}
We split up $|F_{P}(1/2)|^{2}$ here so that the combined contribution from primes $\leq P^{1/C^d}$ inside the large bracket comes with the desired exponent $2(q-1) + 2 = 2q$. As discussed above, for given $d$ we expect $H_{P,d}(\frac{1}{2})$ to be roughly independent of the contribution from primes $> P^{1/C^d}$ inside the large bracket, in other words we expect that $\E H_{P,d}(\frac{1}{2})$ {\em ought} to essentially factor out. This suggests it {\em might} be reasonable to split up $H_{P,d}(\frac{1}{2})$ as $H_{P,d}(\frac{1}{2})^{2-q} H_{P,d}(\frac{1}{2})^{q-1}$, so when we apply H\"older's inequality with exponents $1/(2-q)$ and $1/(q-1)$, we again end up with $H_{P,d}(\frac{1}{2})$ outside the bracket\footnote{This is another place where our heuristics certainly do not guarantee success in advance, and some trial and error is required--- the entire process of applying H\"older's inequality might have turned out to be too wasteful, regardless of the way we split things. Indeed, the larger the power of $H_{P,d}(\frac{1}{2})$ that we move into the bracket the more we gain from the decorrelation between this and $|F_{P}(1/2+ \frac{im}{\log P})|^{2}$, but also the more we risk losing by inbalancing the powers away from $2q$. As it turns out, we obtain bounds where the power of $C^d$ does not quite match the guess we might make about the real size of $\E |F_{P}(1/2)|^{2} \left( \sum_{\substack{C^{d-1} \leq |m| \leq C^d, \\ |m| \leq (\log P)/2}} |F_{P}(1/2+ \frac{im}{\log P})|^{2} \right)^{q-1}$, but is still satisfactory.}.

Proceeding exactly as described, H\"older's inequality implies that for each $d$ we have
\begin{eqnarray}
&& \E |F_{P^{1/C^d}}(\frac{1}{2})|^{2q(2-q)} H_{P,d}(\frac{1}{2})^{2-q} \Biggl( \sum_{\substack{C^{d-1} \leq |m| \leq C^d, \\ |m| \leq (\log P)/2}} |F_{P^{1/C^d}}(\frac{1}{2})|^{2(q-1)} H_{P,d}(\frac{1}{2}) |F_{P}(\frac{1}{2}+ \frac{im}{\log P})|^{2} \Biggr)^{q-1} \nonumber \\
& \leq & \left( \E |F_{P^{1/C^d}}(\frac{1}{2})|^{2q} H_{P,d}(\frac{1}{2}) \right)^{2-q}   \Biggl( \E \sum_{\substack{C^{d-1} \leq |m| \leq C^d, \\ |m| \leq (\log P)/2}} |F_{P^{1/C^d}}(\frac{1}{2})|^{2(q-1)} H_{P,d}(\frac{1}{2}) |F_{P}(\frac{1}{2}+ \frac{im}{\log P})|^{2} \Biggr)^{q-1} . \nonumber
\end{eqnarray}
We can now calculate the expectations, and find this is all
$$ \approx \left( (\frac{\log P}{C^d})^{q^2} C^d \right)^{2-q} \Biggl( \sum_{\substack{C^{d-1} \leq |m| \leq C^d, \\ |m| \leq (\log P)/2}} (\frac{\log P}{C^d})^{q^2} C^{2d} \Biggr)^{q-1} \ll \frac{\log^{q^2}P}{C^{d(q-1)^2}} = \frac{\log^{q^2}P}{e^{d(q-1)}} . $$
Summing over $d$, we get a bound $\ll \frac{\log^{q^2}P}{q-1}$.

We see that this argument delivers the desired bound $\ll \log^{q^2}P$ for $q$ bounded strictly away from 1, but not if $q$ may be close to 1 in a way depending on $P$ (which is possible in Theorem \ref{thmsthigh}). We cannot choose $C$ larger to fix this problem, because when calculating rigorously one loses factors $C^{O(1)}$ in the expectations inside the large bracket (reflecting the fact that the contributions from primes $> P^{1/C^d}$ are not perfectly uncorrelated for the full range of $C^{d-1} \leq |m| \leq C^d$), so we must have $C^{O(q-1)} \ll 1$ when raising that bracket to the power $q-1$.

Instead, we can try to mitigate the inefficiency that enters the argument in \eqref{ddecomp}, by collecting some of the $d$ values together and pulling out a maximum over $d$, rather than a full sum. Working with such a maximum requires a use of {\em martingale theory} and appropriate {\em maximal inequalities}, applied to the filtration structure where one adds batches of successive primes into the partial Euler products $F$. For more information about this, we refer to the original paper~\cite{harperrmfhigh}.

\section{Further reading}
We end by providing a few further references. The study of random multiplicative functions is currently very active: for example, there has been extensive work on distributional results, almost sure bounds, (non-)vanishing results, etc. for sums of $f(n)$, including weighted sums and sums over various interesting subsets of $\N$. Various model settings have been explored as well, including computational work. Beyond this, there is a large body of work developing the connections between random multiplicative functions and number theoretic issues like Dirichlet character sums, moments of $L$-functions, and the Fyodorov--Hiary--Keating conjecture; and probabilistic issues like multiplicative chaos, and secular coefficients of random matrices. A full survey would be far beyond the scope of this paper, so we limit ourselves to pointing out some works rather closely related to Theorems \ref{thmstlow} and \ref{thmsthigh}. The introductions of Harper~\cite{harperrmflowmoments}, of Soundararajan and Zaman~\cite{soundzam}, and of Garban and Vargas~\cite{garbanvargas} provide more detailed overviews (although in some cases now a bit out of date) of some other work in these areas.

A recent paper of Xu~\cite{xubts} adapts the methods underlying Theorem \ref{thmstlow} to study the sum of a random multiplicative function over $R$-rough numbers, finding the threshold for $R$ (in terms of $x$) at which better than squareroot cancellation in the first absolute moment breaks down. Caich~\cite{caichshort} performs a similar investigation for the sum of $f(n)$ over short intervals, finding the threshold of interval length where better than squareroot cancellation appears (and proving a partial analogue for character sums as well).

The weighted sums $\sum_{n \leq x} \frac{f(n)}{n^{\sigma}}$, $\sum_{n \leq x} \frac{f(n) d_{\alpha}(n)}{n^{\sigma}}$ (where $d_{\alpha}$ is a generalised divisor function), and especially $\sum_{n \leq x} \frac{f(n)}{\sqrt{n}}$, have been investigated as possible models for (powers of) the Riemann zeta function. Gerspach~\cite{gerspachpseudo}, and then Gerspach and Lamzouri~\cite{gerspachlam}, adapt some of the methods underlying Theorems \ref{thmstlow} and \ref{thmsthigh} to estimate $\E|\sum_{n \leq x} \frac{f(n) d_{\alpha}(n)}{\sqrt{n}}|^{2q}$ up to lower order factors, on all ranges of $q$ and $\alpha$ where such bounds were not previously known. Aymone, Heap and Zhao~\cite{ahz} refined some of their estimates, amongst various other results exploring different aspects of the behaviour of $\sum_{n \leq x} \frac{f(n)}{\sqrt{n}}$. See also the work of Brevig and Heap~\cite{brevigheap}, who (by different methods) investigate the dependence on $q$ of the implicit constants in estimates for $\E|\sum_{n \leq x} \frac{f(n)}{\sqrt{n}}|^{2q}$ with $q$ large.

Gu and Zhang~\cite{guzhang} have adapted the methods underlying Theorem \ref{thmstlow} to bound the low moments of so-called secular coefficients, in a model sometimes described as holomorphic multiplicative chaos. The Gaussian version of holomorphic multiplicative chaos corresponds to the model setting studied by Soundararajan and Zaman~\cite{soundzam}, and also explored extensively (with a more probabilistic slant) by Najnudel, Paquette and Simm~\cite{najpaqsimm}. Gu and Zhang study non-Gaussian variants of the model, obtaining the same moment bounds as in the Gaussian case (which are analogous to Theorem \ref{thmstlow}) provided their underlying non-Gaussian random variables are sufficiently light tailed. They also prove results, with different behaviour, in certain heavy tailed cases.

Finally, we mention recent work proving analogues of Theorems \ref{thmstlow} and \ref{thmsthigh} for number theoretic averages of Dirichlet character sums and of zeta sums (i.e. of sums $\sum_{n \leq x} \chi(n)$ and $\sum_{n \leq x} n^{it}$, averaging over $\chi$ and over $t$ respectively). Harper~\cite{harpertypicalchar} proves upper bounds for low moments of character and zeta sums. Szab\'o~\cite{szaboupper, szabolower} proves upper and lower bounds for the $2q$-th moment of character sums, when $q > 2$ (and assuming the truth of the Generalised Riemann Hypothesis for the upper bounds). Gao~\cite{gaozeta} recently established a partial analogue for zeta sums of Szab\'o's upper bound result, whilst Baier and Gao~\cite{baiergao} handle upper bounds for character sum moments over function fields. Assuming certain (very strong) conjectures from number theory, Wang and Xu~\cite{wangxu} prove a conjecture of Harper~\cite{harpertypicalchar} on upper bounds for low moments of character sums $\sum_{n \leq x} \chi(n) \lambda(n)$, twisted by the classical Liouville function $\lambda(n)$. This has an application to the distribution of the Liouville function in arithmetic progressions.

\vspace{12pt}
\noindent {\em Acknowledgements.} The author would like to thank Marco Aymone, Ofir Gorodetsky, Mo Dick Wong, Max Xu, and Asif Zaman for their comments and encouragement.

\end{document}